\documentclass[reqno]{amsart}

\usepackage{amsmath,amssymb,amsthm,amsfonts,mathtools}
\usepackage{hyperref}
\usepackage[english]{babel}

\usepackage{cancel}
\usepackage{color}
\usepackage{inputenc}
\usepackage{fontenc}
\usepackage{tikz}

\newcommand{\cD}{{\mathcal D}}

\newcommand{\cH}{{\mathcal H}}
\newcommand{\cL}{{\mathcal L}}

\newcommand{\cQ}{{\mathcal Q}}

\newcommand{\cU}{{\mathcal U}}

\DeclareMathOperator{\wlim}{w-lim}

\newtheorem{theorem}{Theorem}[section]
\newtheorem{lemma}[theorem]{Lemma}
\newtheorem{proposition}[theorem]{Proposition}
\newtheorem{remark}[theorem]{Remark}
\newtheorem{remarks}[theorem]{Remarks}

\newtheorem{assumptions}[theorem]{Assumptions}

\newtheorem{claim}{Claim}
\newtheorem{problem}[theorem]{Problem}
\numberwithin{equation}{section}

\date{}

\begin{document}

\title[]{On the infinitesimal generator of an optimal state semigroup}

\author{Paolo Acquistapace}
\address{Paolo Acquistapace, Universit\`a di Pisa, Dipartimento di Matematica,
Largo Bruno~Pontecorvo 5, 56127 Pisa, ITALY 
}
\email{paolo.acquistapace(at)unipi.it}

\author{Francesca Bucci}
\address{Francesca Bucci, Universit\`a degli Studi di Firenze,
Dipartimento di Matematica e Informatica,
Via S.~Marta 3, 50139 Firenze, ITALY
}
\email{francesca.bucci(at)unifi.it}



\begin{abstract}
In this article we fully describe the domain of the infinitesimal generator of the optimal state semigroup
which arises in the theory of the linear-quadratic problem for a specific class of boundary control systems.
This represents an improvement over earlier work of the authors, joint with Lasiecka, 
where a set inclusion was established, but not an equality. 
The novel part of the proof of this result developes through appropriate asymptotic estimates that take advantage of the regularity analysis carried out in the study of the optimization problem, while 
the powers of positive operators and interpolation are still key tools.
We also attest to the validity of an assumed relation between two significant parameters
in the case of distinct systems of coupled hyperbolic-parabolic partial differential equations which
are pertinent to the underlying framework.
\end{abstract}

\maketitle


\section{Introduction} \label{s:intro}
The linear-quadratic problem is among the most renowned as well as successful research topics within
control theory.
We recall briefly that by infinite horizon linear-quadratic (LQ) problem it is meant the search of
a control function $\hat{u}(\cdot)$ minimizing a quadratic functional such as
\begin{equation} \label{e:cost}
J(u):=\int_0^\infty \left(\|Ry(t)\|_Z^2 + \|u(t)\|_U^2\right)dt\,, 
\end{equation}
within a given class $\cU$ of admissible controls, where $y(\cdot)=y(\cdot;y_0,u)$
is the solution to an initial-value problem for the linear (control) system $y'=Ay+Bu$ on
$(0,\infty)$, that corresponds to $u(\cdot)\in \cU$ and to an arbitrary initial state $y_0\in Y$.
Assuming that such optimal control $\hat{u}(\cdot)$ does exist, then 
$$\hat{y}(\cdot)=y(\cdot;y_0,\hat{u})$$ is termed {\em optimal state}, and the {\sl ensemble} 
$(\hat{y},\hat{u})$ optimal pair.
To frame the discussion that follows, we initially omit the functional-analytic details
pertaining to the spaces $Y$, $U$, $Z$ and the linear operators $A$, $B$ and $R$;
the reader may think of the latter as they all were bounded operators 
(acting between respective Hilbert spaces), just for the sake of simplicity.
However, it should be kept in mind that $A\colon \cD(A)\subset Y \longrightarrow Y$ 
will be always the generator of a $C_0$-semigroup $e^{At}$ in $Y$, $t\ge 0$;
in addition it is assumed that $e^{At}$ is uniformly stable (see the basic Assumptions~\ref{a:ipo_0}). 
And also, it will be $\cU=L^2(0,\infty;U)$.

As it will be made apparent later on, the optimal state satisfies a {\em semigroup} property; to wit,
$\hat{y}(t)=\Phi(t)y_0$, where $\Phi(t)$ is a $C_0$-semigroup on $Y$, $t\ge 0$.
With the focus on a class of abstract control systems which has proven effective to describe systems of
coupled hyperbolic-parabolic partial differential equations (PDE), subject to a boundary control action
on its parabolic component, our goal in this paper is to provide a complete description
of the infinitesimal generator of the said semigroup.
This by means of a full characterization of its domain, 
thus answering to a question that had remained open in our previous works \cite{abl_2013} 
(joint with Lasiecka) and \cite{ab-uniqueness_2020} ({\sl cf.}~\cite[Section~1.1]{ab-uniqueness_2020}).

As it is well known, once the existence of a unique optimal control 
minimizing the cost functional is ascertained -- which may follow easily by using classical variational arguments --, 
the property which is sought and identified as the {\em actual solution} of the LQ problem
is a (pointwise in time) {\em feedback} representation of the optimal control $\hat{u}(\cdot)$ in terms of
the corresponding optimal state $\hat{y}(\cdot)$.
This turns out to be
\begin{equation*}
\hat{u}(t)=-B^*P\hat{y}(t) \qquad \textrm{a.e. in $[0,\infty)$,} 
\end{equation*} 
where the operator $P$ solves the quadratic algebraic equation 
\begin{equation*}
PA+A^*P-PBB^*P+R^*R=0\,, 
\end{equation*}
known as Riccati equation (to be interpreted appropriately, when $B$ is an {\em unbounded} operator; see
\eqref{e:ARE}).
The above Riccati equation is well-posed in a finite dimensional context (where $A$, $B$, $R$ are matrices), and also in the case the abstract state equation describes a boundary value problem for a partial differential equation in a bounded domain with smooth boundary, in the presence of distributed control (which 
naturally yields a {\em bounded} control operator $B$, i.e. such that $B\in \cL(U,Y)$).
A historical synopsis on Riccati equations, along with relevant bibliographical information,
is offered in \cite[Section~1.2]{ab-uniqueness_2020}.

The well-posedness of the Riccati equation is key for the (so-called) synthesis of the optimal control.
Indeed, the latter is accomplished starting from the Riccati equation corresponding to the optimal control problem and singling out its unique solution $P$; next, solving the closed-loop equation 
\begin{equation*}
\begin{cases}
y'=(A-BB^*P)y\,, & \quad t> 0
\\[2mm]
y(0)=y_0\in Y & 
\end{cases}
\end{equation*}
that is obtained taking into account the above feedback form of $\hat{u}(\cdot)$.
Its solution $\hat{y}(\cdot)$ allows to finally determine the optimal control $\hat{u}(\cdot)$.

It was discovered during the eighties of the last century -- thanks to the work of Da Prato and Ichikawa, Flandoli, Lasiecka and Triggiani 
-- that the above process can be extended to control systems which yield an (intrinsically) {\em unbounded} $B$, and yet satisfy the following assumptions.
These are characteristic of parabolic (and parabolic-like) PDE with boundary or point control;
see \cite{bddm}, \cite{las-trig-redbooks}.

\smallskip
{\bf Parabolic class.}
{\em
Let $Y$, $U$ be separable complex Hilbert spaces. 
\begin{itemize}

\item
The closed linear operator $A\colon \cD(A)\subset Y \to Y$ is the infinitesimal generator of a strongly
continous semigroup $\{e^{At}\}_{t\ge 0}$ on $Y$, of type $\omega_0$, which in addition is analytic;

\smallskip

\item 
$B\in \cL(U,[\cD(A^*)]')$; 
there exists $\lambda_0>\omega_0$ and $\gamma\in (0,1)$ such that 
$$(\lambda_0-A)^{-\gamma}B\in \cL(U,Y)\,.$$
\end{itemize}
}
We note that if the semigroup $e^{At}$ underlying the free dynamics is uniformly (exponentially) stable, 
then in particular we may take $\lambda_0=0$.
We highlight several significant features of the (by now, classical) theory of the
LQ problem and of Riccati equations devised within the framework defined by the above hypotheses:
 
\begin{enumerate}

\item[i)]
the optimal state $\hat{y}(t)$ satisfies a semigroup property, that is, 
$\hat{y}(t)=\Phi(t)y_0$, where $\Phi(t)$ is a $C_0$-semigroup on $Y$, $t\ge 0$ 
(which is exponentially stable as well, and in addition {\em analytic} here);

\item[ii)]
there exists a self-adjoint, non-negative definite operator $P\in \cL(Y)$ which satisfies 
$B^*P \in \cL(Y,U)$ and solves the algebraic Riccati equation, that is 
\eqref{e:ARE} below;

\item[iii)]
the Riccati operator $P$ satisfies $(-A^*)^{1-\delta}P\in \cL(Y)$, for any $\delta\in (0,1)$,
a property which is central to the proof of ii);

\item[iv)]
the generator of $\Phi(t)$ is the operator $A_P=A-BB^*P$, whose domain is
\begin{equation*}
\cD(A_P)= \big\{x\in \cD((-A)^{1-\gamma})\colon \,
(-A)^{1-\gamma}x-(-A)^{-\gamma}BB^*Px\in \cD((-A)^{\gamma})\big\}\,;
\end{equation*}
in particular, $\cD(A_P)\subset \cD((-A)^{1-\gamma})$.
\end{enumerate}

\noindent
(See the monographs \cite{bddm} and \cite{las-trig-redbooks}, the latter providing an in-depth treatment.)

\smallskip
When it comes to PDE of hyperbolic type things change dramatically in the worse.
Indeed, the issue which is the most delicate from a theoretical viewpoint, and may be out of reach, 
is that of obtaining a gain operator $B^*P$ which is well-defined or at least densely defined on $Y$,
besides than on the optimal trajectory (as the feedback formula requires by itself).  
See, e.g., the discussion in \cite[Section~1.2.4]{ab-uniqueness_2020} and the cited references.

The case of coupled systems of hyperbolic-parabolic PDE sets itself at the boundary of the two classical
(parabolic and hyperbolic, respectively) classes.
Because significant physical interactions are described by evolutionary PDE systems which comprise dynamics
of different nature -- 
notably, of hyperbolic and parabolic type --, the study of the LQ problem and of the corresponding Riccati equations for these composite PDE systems has received quite a bit of attention over the last twenty years. 
This both at a functional-analytic level, and in the endeavours to pursue 
a clever regularity analysis on a certain interconnected PDE problem (eventually, on a variety of them).
See \cite[Section~1.2]{ab-uniqueness_2020} and the pertinent references.

An abstract class of control systems broad enough to encompass a diverse range of physical interactions has been devised by the authors, jointly with Lasiecka, under well-structured assumptions on 
(the adjoint of) the kernel $e^{At}B$; 
these are detailed in the next subsection as Assumptions~\ref{a:ipo_2}.
A theory of the (finite and) infinite time horizon LQ problem, along with a complex of results
that combine to bring about well-posedness of the corresponding (differential and) algebraic Riccati equations -- thereby attaining the actual synthesis of the optimal control --, has been developed in \cite{abl_2013} and the recent \cite{ab-uniqueness_2020}, following the former study \cite{abl_2005}.
A distinguishing feature of the Riccati theory developed in \cite{abl_2005,abl_2013} is that the gain
operator is bounded on $\cD((-A)^\epsilon)$ (with a suitable $\epsilon>0$), 
which is dense in the state space $Y$;
see the statement A4.~of Theorem~\ref{t:theory-tinfinite}.
Furthermore, it should be noted that the algebraic Riccati equation \eqref{e:ARE}, with $x,y\in \cD(A)$, 
actually extends to be meaningful with $x,y\in \cD(A_P)$; 
see the statement A7. of Theorem~\ref{t:theory-tinfinite}, based on 
\cite[Lemma~5.15 and Theorem~5.16]{abl_2013}.

An element which had remained unresolved in \cite{abl_2013} and \cite{ab-uniqueness_2020} 
was a {\em full} characterization of the domain of the optimal state semigroup's generator $A_P$.
Indeed the information achieved in \cite{abl_2013} with regard to $\cD(A_P)$ was the refined set inclusion 
which follows combining the statement S5. in \cite[Theorem~1.5]{abl_2013} 
(recalled as \eqref{e:inclusion_0} later on) with the property (asserted in \cite[Proposition~5.5]{abl_2013})
$\cD(A_P)\subseteq \cD((-A)^\vartheta)$, valid for any $\vartheta \in (0,1-\gamma)$.
This is
\begin{equation*} 
\cD(A_P)\subseteq E:= \bigcap_{\theta \in (0,1-\gamma)} E_\vartheta\,,
\end{equation*}
having set
\begin{equation} \label{e:set-E_theta}
E_\vartheta:= \big\{x\in \cD((-A)^\vartheta)\colon \, x-A^{-1}B\,B^*Px\in \cD(A) \big\}\,.
\end{equation}
The following question then arises:

\smallskip
{\bf Question.}
{\em
Does the converse $\cD(A_P)\supseteq E$ hold? 
}

\smallskip
\noindent
Our main result, i.e. Theorem~\ref{t:main}, shows that this is indeed the case, 
provided that the parameters $\gamma$ and $\epsilon$ which occur in the Assumptions~\ref{a:ipo_2}
fulfil $\epsilon<1-\gamma$. 
This relation is feasible: to support the assertion, we revisit the actual values of
$\gamma$ and $\epsilon$ brought about by the (model-specific) trace regularity results established in \cite{bucci-las-thermo_2004}, \cite{abl-thermo_2005}, \cite{bucci_2008}, \cite{bucci-dampedfsi_2021}
in regard to distinct PDE problems, in the process of establishing that these fall into
the present underlying framework.

\medskip
The structure of the paper is as follows: in the next subsection we recall the class of control systems
under consideration, together with the outcome of our prior work in \cite{abl_2013} and 
\cite{ab-uniqueness_2020}, summarized as Theorem~\ref{t:theory-tinfinite}.
The statements of Theorem~\ref{t:main} (that is our main result) and of the preliminary
Lemma~\ref{l:inclusion} -- which clarifies a first set inclusion -- are highlighted in the separate Section~\ref{ss:statements}. 

Section~\ref{s:proofs} is mainly devoted to the proof of Theorem~\ref{t:main}.
A neat proof of Lemma~\ref{l:inclusion}, which was not given explicitly in \cite{abl_2013} and 
in addition settles a confusing abuse of notation occurred in \cite[Proposition~5.5]{abl_2013}, is also provided.
  
Finally, in the last section we review three coupled systems of hyperbolic-parabolic PDE describing just as many physical interactions, which attest to the validity of the constraint on the significant parameters
$\gamma$ and $\epsilon$ (presupposed in Theorem~\ref{t:main}).


\subsection{Underlying framework: genesis of the optimal state semigroup}
\label{ss:framework}
Let $Y$ and $U$ be two separable Hilbert spaces, {\em viz.} the {\em state} and {\em control} spaces, respectively.
We consider the infinite-dimensional (linear) control system $y'=Ay+Bu$ on the half-line $[0,\infty)$.
The operator $A$ (which describes the {\em free} dynamics) and the control operator $B$ are initially characterized by the following basic hypotheses.


\begin{assumptions}[Basic hypotheses on the control system] \label{a:ipo_0} 
Let $Y$, $U$ be separable complex Hilbert spaces. 
\begin{itemize}
\item
The closed linear operator $A\colon \cD(A)\subset Y \to Y$ is the infinitesimal generator of a strongly continuous semigroup $\{e^{At}\}_{t\ge 0}$ on $Y$, which is exponentially stable,
i.e. there exist constants $M\ge 1$ and $\omega>0$ such that 
\begin{equation*} 
\|e^{At}\|_{\cL(Y)} \le M \,e^{-\omega t} \qquad \forall t\ge 0\,.
\end{equation*}

\item 
$B\in \cL(U,[\cD(A^*)]')$; equivalently, $A^{-1}B\in \cL(U,Y)$.

\end{itemize}

\end{assumptions}

Then for any given $y_0\in Y$ and any control function $u\in L^2(0,\infty;U)$ the Cauchy problem


\begin{equation}\label{e:cauchy-pbm}
\begin{cases}
y'(t)=Ay(t)+Bu(t)\,, & \quad t> 0
\\[2mm]
y(0)=y_0\in Y & 
\end{cases}
\end{equation}
possesses a unique {\em mild} solution given by
\begin{equation} \label{e:mild}
y(t)= e^{At}y_0+\int_0^t e^{A(t-s)}Bu(s)\,ds\,, \qquad t\ge 0\,.
\end{equation}
The formula \eqref{e:mild} makes sense in the extrapolation space $[\cD(A^*)]'$: 
it yields more specifically $y\in L^2(0,\infty;[\cD(A^*)]')$; 
see \cite[\S\,0.3, p.~6, and Remark~7.1.2, p.~646]{las-trig-redbooks}.

\begin{remark}
\begin{rm}
We note that owing to the first of the Assumptions~\ref{a:ipo_0}, the operator $-A$ is a {\em positive} operator according to \cite[Definition~4.1]{lunardi-interpol}.
This allows to define the powers with complex exponent $(-A)^\alpha$ (in a first step for 
$\text{Re}\,\alpha <0$ via a Dunford integral and then for $\text{Re}\,\alpha \ge 0$); 
see, e.g., \cite[\S~4.1]{lunardi-interpol}.
The fractional powers $(-A)^\vartheta$, $\vartheta\in (0,1)$, will be critically (and repeatedly)
used in the proofs of our results.
We shall write $A^\vartheta$ in place of $(-A)^\vartheta$ throughout the paper 
in order to make the notation lighter.
\end{rm}
\end{remark}


To the state equation in \eqref{e:cauchy-pbm} we associate the quadratic functional 
$J(u)$ defined in \eqref{e:cost},
where the so called {\em observations} space $Z$ is a third separable Hilbert space
(possibly, $Z\equiv Y$). 
It is assumed that $R\in \cL(Y,Z)$ throughout.

\smallskip 
By optimal control problem \eqref{e:cauchy-pbm}-\eqref{e:cost} it is meant the following.


\begin{problem}[\bf Infinite horizon optimal control problem] \label{p:lq-pbm}
Given $y_0\in Y$, seek a control function $u\in L^2(0,\infty;U)$ which minimizes the
cost functional \eqref{e:cost}, where $y(\cdot)=y(\cdot\,;y_0,u)$ is the solution to
\eqref{e:cauchy-pbm} corresponding to the control function $u(\cdot)$ (and with initial state $y_0$)
given by \eqref{e:mild}.
\end{problem}

With motivation coming from the optimal boundary control of partial differential equations systems
which comprise both hyperbolic and parabolic dynamics -- such as, e.g., the ones which describe certain
mechanical-thermal, acoustic-structure, fluid-elasticity interactions -- a thorough study of Problem~\ref{p:lq-pbm} and its complete solution has been provided by the authors in \cite{abl_2013} (jointly with Lasiecka), and in \cite{ab-uniqueness_2020} 
(achieving uniqueness for both differential and algebraic Riccati equations), under the following distinguishing assumptions on the operators $A$, $B$, $R$.
These assumptions were singled out for the most part in \cite{abl_2005}, focused on the finite time horizon optimal control problem, and subsequently strengthened in \cite{abl_2013} to deal with the infinite time horizon one.

\begin{assumptions} \label{a:ipo_2} 
Let $Y$, $U$ and $Z$ be the state, control and observation spaces introduced in the 
Assumptions~\ref{a:ipo_0} and in \eqref{e:cost}, respectively.

The operator $B^*e^{A^*t}$ can be decomposed as 
\begin{equation} \label{e:key-hypo} 
B^*e^{A^*t}x = F(t)x + G(t)x\,, \qquad t\ge 0\,, \; x\in \cD(A^*)\,,
\end{equation}
where $F(t)\colon Y\longrightarrow U$ and $G(t)\colon \cD(A^*)\longrightarrow U$, $t>0$, are bounded linear
operators satisfying the following assumptions:

\begin{enumerate}
\item[1.] 
there exist constants $\gamma\in (0,1)$ and $N, \eta>0$ such that 
\begin{equation} \label{e:singular}
\|F(t)\|_{\cL(Y,U)} \le N\,t^{-\gamma}\,e^{-\eta t}\qquad \forall t>0\,,
\end{equation}

\item[2.] 
there exists $T>0$ such that the operator $G(\cdot)$ belongs to $\cL(Y,L^p(0,T;U))$ for all $p\in [1,\infty)$;
\item[3.] 
with $T$ as above, there exists $\epsilon>0$ such that:
\begin{enumerate}
\item[a)] 
the operator $G(\cdot){A^*}^{-\epsilon}$ belongs to $\cL(Y,C([0,T];U))$, 
with
\begin{equation*} 
\sup_{t\in [0,T]}\|G(t){A^*}^{-\epsilon}\|_{\cL(Y,U)} <\infty\,;
\end{equation*}
\item[b)] 
the operator $R^*R$ belongs to $\cL(\cD(A^{\epsilon}),\cD({A^*}^{\epsilon}))$, i.e.
\begin{equation*} 
\|{A^*}^{\epsilon}R^*RA^{-\epsilon}\|_{\cL(Y)} \le c<\infty\,;
\end{equation*}
\item[c)] 
there exists $q\in (1,2)$ (depending, in general, on $\epsilon$) such that the map 
$x \longmapsto B^*e^{A^*t}{A^*}^\epsilon x$
has an extension which belongs to $\cL(Y,L^q(0,T;U))$. 
\end{enumerate}

\end{enumerate}

\end{assumptions}

\smallskip
As already pointed out in the Introduction, 
the asserted {\em complete solution} to Problem~\ref{p:lq-pbm} embodies several major achievements.
We recall that the existence of a unique optimal pair $(\hat{u},\hat{y})$ follows readily from convex optimization arguments. 
The soughtafter {\em closed loop} form of the unique optimal control $\hat{u}$, i.e. \eqref{e:feedback} below, holds true, as well. 
However, its proof is far from being straightforward: it requires further analytical work than the 
classical argument based on the optimality conditions.
Indeed, the tricky bit has been to ascertain that given the optimal cost operator 
$P\in \cL(Y)$ -- intrinsically defined in terms of the optimal state $\hat{y}(t) = \Phi(t)y_0$ (just like in previous theories, see \eqref{e:optimal-cost-op} below) --, then the {\em gain operator} $B^*P$ is meaningful and bounded at least on a suitable dense subset of the state space $Y$ (which turns out to be $\cD(A^\epsilon)$), and then
to prove that $P$ does solve the algebraic Riccati equation corresponding to Problem~\ref{p:lq-pbm}, displayed below as \eqref{e:ARE}.
All this 
\begin{itemize}
\item
in the absence of the smoothing properties inherited by both the optimal evolution
and the Riccati operator in the fully parabolic case (highlighted in Section~\ref{s:intro} inside i) and 
as iii), respectively), and also 
\item
under a very minimal additional assumption on the observation operator $R$ (that is 3b) of the Assumptions~\ref{a:ipo_2}).
\end{itemize}
In the following Theorem we collect all the pertinent assertions;
we refer the reader to the part titled {\sl Present work} in the Introduction of \cite{abl_2013}, for a synopsis of the analysis carried out therein, as the paper is long and technical.

\begin{theorem}[\cite{abl_2013}, Theorem~1.5; \cite{ab-uniqueness_2020}, Theorem~2.11] 
\label{t:theory-tinfinite}
Under the Assumptions~\ref{a:ipo_2}, the following statements are valid.
\begin{enumerate}

\item[\bf A1.] 
For any $y_0\in Y$ there exists a unique optimal pair 
$(\hat{u}(\cdot),\hat{y}(\cdot))$ for Problem~\eqref{e:cauchy-pbm}-\eqref{e:cost}, 
which satisfies the following regularity properties
\begin{align*}
& \hat{u}\in \bigcap_{2\le p<\infty} L^p(0,\infty;U)\,,
\\
& \hat{y}\in C_b([0,\infty);Y) \cap \Big[\bigcap_{2\le p< \infty} L^p(0,\infty;Y)\Big]\,.
\end{align*}

\item[\bf A2.] 
The family of operators $\Phi(t)$, $t\ge 0$, defined by 
\begin{equation} \label{e:optimal-state-semigroup}
\Phi(t)y_0 :=\hat{y}(t)=y(t,y_0;\hat{u})
\end{equation}
is a $C_0$-semigroup on $Y$, $t\ge 0$, which is exponentially stable.

\item[\bf A3.] 
The operator $P\in\cL(Y)$ defined by 
\begin{equation} \label{e:optimal-cost-op}
Py_0 := \int_0^\infty e^{A^*t}R^*R \Phi(t)y_0\,dt\,, \qquad x\in Y\,,
\end{equation}
is the {\em optimal cost operator}; $P$ is (self-adjoint and) non-negative.

\item[\bf A4.] 
The following ({\em pointwise in time}) feedback representation of the optimal control 
is valid for any initial state $y_0\in Y$:
\begin{equation} \label{e:feedback}
\hat{u}(t) = - B^*P \hat{y}(t) \qquad \textrm{for a.e. $t\in (0,\infty)$},
\end{equation}
where the gain operator satisfies $B^*P\in \cL(\cD(A^\epsilon),U)$
(that is, it is just densely defined on $Y$ and yet it is bounded on $\cD(A^\epsilon)$).

\item[\bf A5.] 
The domain of the infinitesimal generator $A_P$ of the optimal evolution 
$\Phi(t)$ (defined in \eqref{e:optimal-state-semigroup}) satisfies 
\begin{align} 
& \qquad\cD(A_P) \subset \big\{x\in Y\colon  x-A^{-1}B\,B^*Px\in \cD(A) \big\}
\label{e:inclusion_0}\\[1mm]
& \qquad\qquad\quad \subset \Big\{x\in Y\colon 
\exists \wlim_{t\to 0^+} \frac{1}{t}\int_0^t e^{A(t-\tau)}A^{-1}B\,B^*P\Phi(\tau)x\,d\tau
\, \text{in $Y$}\Big\}
\notag
\end{align}
and $A_P\equiv A(I-A^{-1}BB^*P)$ on $\cD(A_P)$.

\item[\bf A6.] 
The operator $e^{At}B$, defined in $U$ and a priori with values in $[\cD(A^*)]'$,
is such that
\begin{equation*} \label{e:tricky-regularity}
e^{\delta\cdot}e^{A\cdot}B\in \cL(U,L^p(0,\infty;[\cD({A^*}^{\epsilon})]') 
\quad 
\forall p\in [1,1/{q})
\end{equation*}
for all $\delta\in [0,\omega\wedge \eta)$; 
almost the very same regularity is inherited by the operator $\Phi(t)B$:
\begin{equation*} \label{e:low-regularity}
e^{\delta\cdot}\Phi(\cdot)B\in \cL(U,L^p(0,\infty;[\cD({A^*}^{\epsilon})]') 
\quad 
\forall p\in [1,1/{q})\,,
\end{equation*}
provided $\delta\in [0,\omega\wedge \eta)$ is sufficiently small.

\item[\bf A7.] 
The optimal cost operator $P$ defined in \eqref{e:optimal-cost-op} is a solution to 
the algebraic Riccati equation (ARE) corresponding to Problem~\eqref{e:cauchy-pbm}-\eqref{e:cost}, 
that is 
\begin{equation} \label{e:ARE}
\begin{split}
& \qquad (Px,Az)_Y+(Ax,Pz)_Y-(B^*Px,B^*Pz)_U+(Rx,Rz)_Z=0\,, 
\\[1mm]
& \hspace{85mm} x,z\in \cD(A)\,;
\end{split}
\end{equation}
the ARE reads as
\begin{equation*}
(A^*Px,z)_Y+(x,A^*Pz)_Y-(B^*Px,B^*Pz)_U+(Rx,Rz)_Z=0  
\end{equation*}
when $x,z\in \cD(A_P)$.

\item[\bf A8.] 
The algebraic Riccati equation \eqref{e:ARE} has a unique solution $P$ within the class
\begin{equation} \label{e:classeQ}
\cQ :=\big\{ Q\in \cL(Y)\colon Q=Q^*\ge 0\,, \; B^*Q\in \cL(\cD(A^\epsilon),U)\big\}\,;
\end{equation} 
hence, the optimal cost operator $P$ defined by \eqref{e:optimal-cost-op} is precisely that solution.

\end{enumerate}
\end{theorem}


\subsection{The generator of the semigroup: statement of the main result}
\label{ss:statements}
The assertion A5.~in Theorem~\ref{t:theory-tinfinite} is the one which specifically pertains to
the generator $A_P$ of the optimal state semigroup $\Phi(t)$.
Although a full characterization of $\cD(A_P)$ was not needed in \cite{abl_2013} (and actually absent from it), the basic set inclusion \eqref{e:inclusion_0} was made more precise therein. 
Indeed, \cite[Proposition~5.5]{abl_2013} related the domain of $A_P$ (and $A_P^*$) to the domain of certain fractional powers of the positive operator $-A$ ($-A^*$, respectively).
For the readers' convenience we explicitly recall the aforesaid result below.

\begin{proposition}[Cf.~\cite{abl_2013}, Proposition~5.5]\label{p:one-inclusion}
With reference to the optimization Problem~\eqref{e:cauchy-pbm}-\eqref{e:cost}, ~
let the Assumptions~\ref{a:ipo_2} be in place.
Accordingly, $\gamma$ is the exponent in \eqref{e:singular}.
Then, the following inclusions are valid, for any $\vartheta\in (0,1- \gamma)$:
\begin{equation} \label{e:link-to-powers}
\cD(A_P) \subset \cD(A^\vartheta)\,, \qquad \cD(A_P^*) \subset \cD({A^*}^\vartheta)\,.
\end{equation}
\end{proposition}

\smallskip
\begin{remarks}
\begin{rm}
i) We note that the letter $\epsilon$ utilized in (5.8) of \cite[Proposition~5.5]{abl_2013}
to indicate the fractional power of $-A$, in place of the chosen $\vartheta$ in \eqref{e:link-to-powers} above,
was certainly misleading.
This is because $\epsilon$ has a clear, specific meaning in 3a-3b-3c. of the Assumptions~\ref{a:ipo_2},
and then in the statements A4. and A8. of Theorem~\ref{t:theory-tinfinite} established under those
hypotheses.

\noindent
ii) It is important to emphasize that the inclusion in the right of \eqref{e:link-to-powers} was 
crucially employed in the proof of \cite[Corollary 5.14]{abl_2013}, where we showed that $\Phi(t)=e^{A_P t}$ is differentiable on $\cD(A)$ 
(a result which is not expected, since $\cD(A)$ is not a natural domain of the strongly perturbed 
evolution $\Phi(t)$).
This regularity property, in turn, was the key to proving that the optimal cost operator $P$ does satisfy
the algebraic Riccati equation; see \cite[Theorem~5.16]{abl_2013}.  
\end{rm}
\end{remarks}

Having recalled from \cite{abl_2013} the findings stated above as Proposition~\ref{p:one-inclusion}, 
and more specifically the first of the inclusions \eqref{e:link-to-powers} therein, 
it becomes apparent that the original information \eqref{e:inclusion_0} on $\cD(A_P)$ provided by Theorem~\ref{t:theory-tinfinite} was already made more precise.
It is worth highlighting this result stemming from \cite{abl_2013}, by means of a separate statement.


\begin{lemma}[Set inclusion]\label{l:inclusion}
With reference to the optimization Problem~\eqref{e:cauchy-pbm}-\eqref{e:cost}, let the Assumptions~\ref{a:ipo_2} be in place; accordingly, $\gamma$ is the exponent in \eqref{e:singular}.
Let $P\in \cL(Y)$ and $A_P\colon \cD(A_P)\subset Y\longrightarrow Y$ denote the corresponding optimal cost operator and the generator of the optimal state semigroup $\{e^{A_Pt}\}_{t\ge 0}$, 
respectively.  
Then,
\begin{equation}\label{e:inclusion_1}
\cD(A_P) \subseteq 
\bigcap_{\vartheta\in (0,1-\gamma)} 
\big\{x\in \cD(A^\vartheta)\colon  x-A^{-1}B\,B^*Px \in \cD(A)\big\}\,.
\end{equation}

\end{lemma}

\smallskip
The question as to whether the set inclusion \eqref{e:inclusion_1} might be enhanced to become an equivalence
was (wrongly) dismissed in the end of \cite[Remark~5.4]{abl_2013}.
The present study reveals that the answer is actually positive, provided $\epsilon$ is 
below a certain threshold determined by $\gamma$.


\begin{theorem}[Main result]\label{t:main}
With reference to the optimization Problem~\eqref{e:cauchy-pbm}-\eqref{e:cost}, let the Assumptions~\ref{a:ipo_2} be in place;
accordingly, $\gamma$ and $\epsilon$ are the exponents 
therein.
Let $P\in \cL(Y)$ and $A_P\colon \cD(A_P)\subset Y\longrightarrow Y$ denote the corresponding optimal cost operator and the generator of the optimal state semigroup $\{e^{A_P t}\}_{t\ge 0}$, 
respectively.  
If $\epsilon <1-\gamma$, then
\begin{equation}\label{e:characterization}
\cD(A_P) \equiv \bigcap_{\vartheta\in (0,1-\gamma)}  
\big\{x\in \cD(A^\vartheta)\colon  x-A^{-1}B\,B^*Px \in \cD(A)\big\}\,.
\end{equation}

\end{theorem}

The proof of Theorem~\ref{t:main}, as well as the one of Lemma~\ref{l:inclusion}, are given in the next section.


\section{From a set inclusion to a set equality}
\label{s:proofs}
This section is devoted to showing our main result.
First we give a sketch of the proof of Lemma~\ref{l:inclusion}, with focus on the first one of the inclusions \eqref{e:link-to-powers}. (We remind the reader that the proof of \cite[Proposition~5.5]{abl_2013} pertained
exclusively to the second one.)
Thus we take the chance of amending a few vague or imprecise assertions contained in the proof of 
\cite[Proposition~5.5]{abl_2013}.
Also, we find it worth emphasizing the powerful role of the theory of interpolation spaces, which
specifically provides the tools. 
Next, we give the proof of Theorem~\ref{t:main}, via a new inclusion (the converse), 
thereby achieving the set equality \eqref{e:characterization}.

\smallskip
It is useful to recall the notation $L$ for the operator -- also termed the ``input-to-state map'' -- defined by
\begin{equation} \label{e:input-to-state}
L\colon u(\cdot) \longrightarrow 
(L u)(\cdot) :=\int_0^{\cdot} e^{A(\cdot-r)}Bu(r)\,dr\,. 
\end{equation}
The assumed decomposition \eqref{e:key-hypo} for the adjoint of the kernel $e^{At}B$ yields
the (consistent) decomposition
\begin{equation}\label{e:splitting}
Lu(t)=\int_0^t F(t-r)^*Bu(r)\,dr+\int_0^t G(t-r)^*u(r)\,dr=: L_{(1)}u(t) +L_{(2)}u(t)\,,
\end{equation}
which will be directly 
useful in the proof of Lemma~\ref{l:inclusion}.
(The notation $L_{(i)}$ for the above operators was introduced in \cite{abl_2013}; see (3.10) therein;
we will write $L_{(i)}u(t)$ in place of $(L_{(i)}u(\cdot))(t)$ for the sake of simplicity.)


\smallskip
\noindent 
{\em Proof of Lemma~\ref{l:inclusion}.} 
The inclusion \eqref{e:inclusion_1} that we are asked to show requires that we compare the domain of the optimal state generator $A_P$ with the domain of a fractional power $A^{\vartheta}$.
Here the key tool is interpolation.
More specifically,

\begin{enumerate}

\item[i)]
we will use the following relation between interpolation spaces and domains of fractional powers:
\begin{equation} \label{e:powers}
(X,\cD(A))_{\vartheta,1}\subseteq \cD(A^{\vartheta}) \subseteq (X,\cD(A))_{\vartheta,\infty}\,, 
\qquad \vartheta\in (0,1)\,,
\end{equation}
whose validity is ensured here by the exponential stability of the semigroup $e^{At}$; 
see the statement (d) in \cite[\S~1.15.2]{triebel}, or \cite[Prop.~4.7, p.~92]{lunardi-interpol}. 
 
\item[ii)]
In turn, the driving idea is to relate $\cD(A_P)$ to an appropriate interpolation space 
$(Y,\cD(A))_{\alpha,\infty}$, first.
This is because the characterization
(\cite[Theorem 1.13.2]{triebel}, \cite[Prop.~5.7]{lunardi-interpol})
\begin{equation} \label{e:interpolation} 
(Y,\cD(A))_{\alpha,\infty}=\Big\{x \in Y \colon 
\sup_{t\in (0,1]} t^{-\alpha}\|e^{At}x-x\|_Y < \infty\Big\}
\end{equation}
is especially suited for the task at hand.

\item[iii)]
The conclusive argument relies on the series of set inclusions
\begin{equation}\label{e:queue}
(X,D(A))_{\alpha,1}\subseteq (X,D(A))_{\alpha,p}\subseteq (X,D(A))_{\alpha,\infty} 
\subseteq (X,D(A))_{\vartheta,1}\,,
\end{equation}
which hold true for any $\alpha, \vartheta, p$ such that $0 <\vartheta < \alpha < 1$ and $1 < p < \infty$;
see \cite[\S~1.3.3]{triebel} or \cite[Propositions~1.3 and 1.4]{lunardi-interpol}.

\end{enumerate}

\smallskip
Let $x\in \cD(A_P)$ be given. 
As pointed out in ii), we will show that $x \in (Y,\cD(A))_{\alpha,\infty}$
for any $\alpha \in (0,1-\gamma)$.
By \eqref{e:mild} with $x$ as initial state, we know that
\begin{equation}\label{e:from-mild} 
(e^{At}-I)x = (\Phi(t)-I)x - (L\hat{u}(\cdot,x))(t)\,, 
\end{equation}
where $\Phi(t)x$ is the optimal state and $\hat{u}(\cdot,x)$ is the optimal control, respectively.
Since the asymptotic estimate $(\Phi(t) - I)x = O(t)$, as $t \to 0^+$, is trivial, 
we focus on the subsequent integral term $L\hat{u}(t)$.
For this we need the splitting \eqref{e:splitting}.
The computations leading to the estimates for $L_{(1)}\hat{u}(t)$ and $L_{(2)}\hat{u}(t)$ are essentially the ones in (5.11) and (5.14) of \cite[Proposition~5.5]{abl_2013}, respectively.
We omit most part of the details and only highlight the principal tools and possible
constraints on the parameters' values.   
 
To produce an estimate of $L_{(1)}\hat{u}(t)$, we recall the (singular) estimate \eqref{e:singular} 
pertaining to the component $F(t)$ of $B^*e^{A^*t}$ (i.e.~1. of the Assumptions~\ref{a:ipo_2}), 
as well as the regularity property
\begin{equation*}
\|\hat{u}\|_{L^p(0,\infty;U)}\le c_p\, \|x\|_Y \qquad \forall p\in [2,\infty)
\end{equation*}
for the optimal control (that is (4.3b) in \cite[Proposition~4.1]{abl_2013}, proved therein).
The above estimate actually extends to any $p\ge 1$ if we limit ourselves to $L^p(0,T;U)$, with $T$ as
in 2.-3. of the Assumptions~\ref{a:ipo_2}.
This makes sense since we are interested in the behaviour of $L_{(i)}\hat{u}(t)$, $i=1,2$, only in
a neighbourhood of $t=0$.    
Then we have
\begin{equation}\label{e:L1-bound}
\|L_{(1)}\hat{u}(t))\|_Y 
\le \int_0^t \frac{C_1}{(t-s)^{\gamma}} \|\hat{u}(s)\|_U\,ds 
\le C_1 t^{1/p'-\gamma} \|x\|_Y\,,
\end{equation}
via the H\"older inequality in the latter bound, where $p'$ is the conjugate exponent of $p$.
We note that in order to have $1/p'-\gamma>0$ in \eqref{e:L1-bound}, one needs $p > 1/(1-\gamma)$.
This is allowed, because $p$ can be taken arbitrarily large; accordingly, $p'$ will be arbitrarily close 
to $1$, which yields the {\em upper threshold} $1-\gamma$ for the powers of $t$.
Therefore, \eqref{e:L1-bound} establishes
\begin{equation}\label{e:L1-final} 
\|L_{(1)}\hat{u}(t)\|_Y\le C_1 t^{\alpha_1}\|x\|_Y
\end{equation}
{\em for any} $\alpha_1 \in (0,1-\gamma)$.

\smallskip
The analysis of the summand $L_{(2)}\hat{u}(t)$ is slightly more delicate, while leading to a better estimate.
We take the inner product with an arbitrary $z\in Y$, use once more that $\hat{u}\in L^p(0,T;U)$ for any finite $p\ge 1$, this time combined with the basic regularity of the component $G(t)$ of $B^*e^{A^*t}$, that is 2a. of the Assumptions~\ref{a:ipo_2}.
An application of the H\"older and Young inequalities gives
\begin{equation}\label{e:L2-bound}
\begin{split}
|(L_{(2)}\hat{u}(t),z)_Y| 
& = \Big|\int_0^t (\hat{u}(s),G(t-s)z)_U\, ds \Big|
\\[1mm]
& \le C_2\, t^{1/r} \|x\|_Y \|z\|_Y\,,\qquad  z\in Y\,, \; t\in [0,T]\,,
\end{split}
\end{equation}
with $r\in (1,\infty)$ that can be taken arbitrarily close to $1$.
Then \eqref{e:L2-bound} establishes
\begin{equation}\label{e:L2-final} 
\|L_{(2)}\hat{u}(t)\|_Y \le C_2\, t^{\alpha_2} \|x\|_Y
\end{equation}
{\em for any} positive $\alpha_2<1$. 

We return to \eqref{e:from-mild} with the information provided by the estimates \eqref{e:L1-final} and \eqref{e:L2-final} to find 
\begin{equation*}
\|(e^{At}-I)x\|_Y = O(t)-O(t^\alpha) = O(t^\alpha)\,, \qquad  t \to 0^+\,,
\end{equation*}
valid for $\alpha :=\alpha_1\wedge \alpha_2$, and hence for any $\alpha\in (0,1-\gamma)$. 
By \eqref{e:interpolation}, we have so far proved the membership 
\begin{equation*}
x\in (Y,\cD(A))_{\alpha,\infty} \qquad \forall \alpha\in (0,1-\gamma)\,.
\end{equation*}
Now we recall \eqref{e:queue} and \eqref{e:powers} to deduce that 
\begin{equation*}
x \in (X,\cD(A))_{\vartheta,1} \subseteq \cD(A^{\vartheta})
\end{equation*}
for any $\vartheta\in (0,\alpha)$ (and any $\alpha\in (0,1-\gamma)$).
Since $x$ was given and yet arbitrary, then $\cD(A_P)\subseteq \cD(A^{\vartheta})$ 
fo any $\vartheta\in (0,1-\gamma)$. 
Finally, the inclusion \eqref{e:inclusion_0} formerly established in \cite{abl_2013} ensures that 
any $x\in \cD(A_P)$ fulfils $x-A^{-1}BB^*Px\in \cD(A)$, which confirms the full inclusion \eqref{e:inclusion_1}, thus concluding the proof of the Lemma.

\qed

\medskip 
The proof of the reverse inclusion requires a different argument, with the definition of infinitesimal
generator of a $C_0$-semigroup as natural starting point (\cite{pazy}).
The proof develops through several distinct steps, with a {\em bootstrap} argument playing a key role
among them. 

 
\smallskip 
\noindent 
{\em Proof of Theorem~\ref{t:main}.} 
In view of Lemma~\ref{l:inclusion}, in order to attain the full characterization \eqref{e:characterization} 
we seek to prove the converse of \eqref{e:inclusion_1}.

\smallskip 
\noindent
{\bf 0.} 
The beginning of the proof is based upon the following observation: 
since $\cD(A^{\vartheta_1})\subseteq \cD(A^{\vartheta_2})$ for $\theta_1\ge \theta_2$, 
then $E_{\vartheta_1}\subseteq E_{\vartheta_2}$ as well (the notation $E_\theta$ has
been introduced in \eqref{e:set-E_theta}).
Therefore, if $\epsilon<1-\gamma$, we have
\begin{equation*}
E:=\bigcap_{\theta \in (0,1-\gamma)} E_\vartheta 
\equiv \bigcap_{\theta \in [\epsilon,1-\gamma)} E_\vartheta=: F\,.
\end{equation*}
Thus, proving 
\begin{equation}\label{e:seconda}
\cD(A_P) \supseteq E_\vartheta \qquad \forall \vartheta \in [\epsilon,1-\gamma)
\end{equation}
will yield $\cD(A_P) \supseteq F\equiv E$, as desired.

Now if $\vartheta\ge \epsilon$ we know that $\cD(A^\vartheta)\subseteq\cD(A^\epsilon)$,
which gives a meaning to the element $A^{-1}BB^*Px$ for $x\in \cD(A^\vartheta)$,
since $B^*P\in \cL(\cD(A^\epsilon),U)$ (by the statement A4. of Theorem~\ref{t:theory-tinfinite}),
while $A^{-1}B\in \cL(U,Y)$ (by the Assumptions~\ref{a:ipo_0}).
 
\smallskip 
\noindent
{\bf 1.} 
Let $\vartheta \in [\epsilon,1-\gamma)$ be given, and let $x\in E_\vartheta$.
By the very definition of $E_\vartheta$, we have $x\in \cD(A^\vartheta)$ and $x-A^{-1}BB^*Px \in \cD(A)$.
To infer that $x\in \cD(A_P)$, we need to address the existence of the limit of $$\frac{\Phi(t)-I}{t}x$$
in $Y$, as $t\to 0^+$, where $\Phi(t)x=\hat{y}(t)$ is the optimal state of the minimization problem.
By the representation \eqref{e:mild} of any mild solution to the Cauchy problem \eqref{e:cauchy-pbm},
with $x$ (in place of $y_0$) as initial state and $\hat{u}(\cdot)$ in its feedback form \eqref{e:feedback}
as the (optimal) control, we know that the optimal state semigroup $\Phi(t)$ satisfies the integral equation
\begin{equation*}
\Phi(t)x =e^{At}x - \int_0^t e^{A(t-s)}BB^*P\Phi(s)x\,ds\,, \qquad t\ge 0\,.
\end{equation*}
Then,
\begin{align} \label{antefatto}
(\Phi(t)-I)x 
&=(e^{At}-I)x - \int_0^t e^{A(t-s)}BB^*P\Phi(s)x\,ds
\notag 
\\
&=\underbrace{(e^{At}-I)(x-A^{-1}BB^*Px)}_{T_1(t,x)} 
+ \underbrace{\int_0^t e^{A(t-s)}BB^*P(\Phi(s)-I)x\,ds}_{T_2(t,x)}
\end{align}
for any $t\ge 0$.
The identity \eqref{antefatto} yields, for $t>0$,
\begin{equation} \label{start}
\frac{\Phi(t)-I}{t}x = \frac{e^{At}-I}{t}(x-A^{-1}BB^*Px)
+ \frac{1}{t}\int_0^t e^{A(t-s)}BB^*P(\Phi(s)-I)x\,ds\,.
\end{equation}
Because $x-A^{-1}BB^*Px\in \cD(A)$, the first summand in the right hand side of \eqref{start}
tends to $A(x-A^{-1}BB^*Px)$, as $t\to 0^+$.
If we can establish that the second summand is such that
\begin{equation*}
\lim_{t\to 0^+}\frac{T_2(t,x)}{t}=0\,,
\end{equation*}
then from \eqref{start} it will follow $x\in \cD(A_P)$, $A_Px=A(x-A^{-1}BB^*Px)$.

The goal of the argument that follows is to show that $T_2(t,x)$ defined in \eqref{antefatto}
is an $\textrm{O}(t^\alpha)$ with $\alpha>1$, as $t \to 0^+$. 

\medskip
\noindent
{\bf 2.}
We begin with some preliminary considerations, which in particular highlight the {\em space} regularity
of the incremental ratio.
Based on \cite[Proposition 4.7]{abl_2013}, from 
$x\in \cD(A^\vartheta)\subseteq \cD(A^\epsilon)$ it follows $\Phi(t)x\in \cD(A^\epsilon)$, which gives
$(\Phi(t)-I)x\in \cD(A^\epsilon)$ as well.

As for the right hand side of \eqref{antefatto}, it is readily seen that the first summand $T_1(t,x)$
belongs to $\cD(A)$, since $x-A^{-1}BB^*Px \in \cD(A)$ (just use a basic property of operator semigroups).
Therefore, clearly the term $T_2(t,x)$ belongs to $\cD(A^\epsilon)$ as well, for any $t\ge 0$.  
However, for the purposes of the estimates we are about to carry out,  
we rewrite $T_2(t,x)$ more cleanly as
\begin{equation*}
T_2(t,x) = L(B^*P(\Phi(\cdot)-I)x)(t)\,,
\end{equation*}
and observe that 
\begin{itemize}

\item
the property $B^*P\in \cL(\cD(A^\epsilon),U)$ combined with the exponential stability of $\Phi(t)$ 
({\em cf.}~the assertions A4. and A2. of Theorem~\ref{t:theory-tinfinite}) ensure in particular that
\begin{equation*} 
B^*P(\Phi(\cdot)-I)x \in C_b([0,\infty),U)\cap L^{q'}(0,\infty;U)\,;
\end{equation*}

\smallskip

\item
$L\in \cL(L^{q'}(0,\infty;U),C_b([0,\infty),\cD(A^\epsilon))$ by \cite[Proposition~3.6(v)]{abl_2013}\,.

\end{itemize}

\noindent
These regularity properties are somewhat excessive. 
Indeed, with respect to the {\em time} regularity, what matters is the behaviour in a right neighbourhood of $t=0$.
 
In the next step we move on to produce (successively enhanced) asymptotic estimates for $(\Phi(t)-I)x$.
Since we are interested in the asymptotic behaviour of the said function (of $t$), as $t\to 0^+$, 
we are allowed to assume $t \le T$, where $T$ is as in 2. of the Assumptions~\ref{a:ipo_2}. 
We will denote by $c$ any positive constant not depending on the variable $t$.

\smallskip
\noindent
{\bf 3.} 
We begin with the derivation of a first estimate for $(\Phi(t)-I)x$ around $t=0$.
The considerations above justify the use of the stronger topology of the space $\cD(A^\epsilon)$.


\begin{claim} \label{claim1} 
The following estimate is valid: 
\begin{equation}\label{base} 
\|(\Phi(t)-I)x\|_{\cD(A^\epsilon)} \le c \big(t^{1-\epsilon}\vee t^{1/q'}\big)\,, \qquad 0\le t\le T\,,
\end{equation}
where $q'$ is the conjugate exponent of $q$ in 3c. of the Assumptions~\ref{a:ipo_2}. 
\end{claim}
To confirm Claim~\ref{claim1}, let us examine either term $T_i(t,x)$, $i=1,2$, 
in the right hand side of \eqref{antefatto}.
We estimate
\begin{equation} \label{stimabase}
\begin{split}
\big\|T_1(t,x)\big\|_{\cD(A^\epsilon)} 
&=\big\|(e^{At}-I)(x-A^{-1}BB^*Px)\big\|_{\cD(A^\epsilon)} 
\\[1mm]
&= \big\|(e^{At}-I)\,A^\epsilon (x-A^{-1}BB^*Px)\big\|_Y\le c\, t^{1-\epsilon}\,, 
\end{split}
\end{equation}
with the bound following by interpolation, since $A^\epsilon (x-A^{-1}BB^*Px)\in \cD(A^{1-\epsilon})$.
Next, in view of $B^*P\in \cL(\cD(A^\epsilon),U)$ (already recalled in the step~2.), we have
\begin{equation*}
B^*P(\Phi(\cdot)-I)x\in C([0,T],U)\subset L^{q'}(0,T;U)\,,
\end{equation*}
so that 
\begin{equation*}
T_2(\cdot,x) =[L(B^*P(\Phi(\cdot)-I)x)](\cdot)\in C([0,T],\cD(A^\epsilon))\,,
\end{equation*}
along with the estimate
\begin{equation} \label{stimabase2}
\begin{split}
\big\|T_2(t,x)\big\|_{\cD(A^\epsilon)}  
&=\left\|A^\epsilon \int_0^t e^{A(t-s)}BB^*P(\Phi(s)-I)x\, ds\right\|_Y 
\\[1mm]
&\le \|(\Phi(\cdot)-I)x\|_{L^{q'}(0,t;\cD(A^\epsilon))}
\\[1mm]
&
\le \|(\Phi(\cdot)-I)x\|_{C([0,t],\cD(A^\epsilon))} t^{1/q'}
\end{split}
\end{equation}
for any $t\in [0,T]$.
Combining \eqref{stimabase2} with \eqref{stimabase} we attain the estimate \eqref{base}
for the left hand side of \eqref{antefatto}, as needed.


\smallskip
\noindent
{\bf 4.}
Let now $n_0=[q'(1-\epsilon)]$ be the largest integer not exceeding $q'(1-\epsilon)$. 
If $n_0=0$, i.e. $q'< \frac1{1-\epsilon}$, then from \eqref{base} it follows
\begin{equation*} 
\|(\Phi(t)-I)x\|_{\cD(A^\epsilon)} \le c\,t^{1-\epsilon}\,, \qquad 0\le t \le T\,.
\end{equation*}
If $n_0>0$ instead, we use a bootstrap argument, based on the following fact.

  
\begin{claim} \label{bootstrap} 
Let $\alpha$ be a positive real number.
If
\begin{equation*}
\left\|A^\epsilon(\Phi(t)-I)x \right\|_Y \le c\, t^{\alpha}, \qquad 0\le t\le T\,,
\end{equation*}
then we have
\begin{equation*}
\left\|A^\epsilon\int_0^t e^{A(t-s)}BB^*P(\Phi(s)-I)x\, ds\right\|_Y \le c\, t^{\alpha + \frac1{q'}}\,, 
\qquad 0\le t\le T\,.
\end{equation*}

\end{claim}
To show that Claim~\ref{bootstrap} holds true, we return to the right hand side of \eqref{antefatto}.
Still with the bound \eqref{stimabase} for $T_1(t,x)$, we aim to improve the estimate 
\eqref{stimabase2} for $T_2(t,x)$.
We argue by duality, using 3c. of the Assumptions~\ref{a:ipo_2}. 
Let $z\in Y$ be given: then, for $0\le t\le T$, we have
\begin{equation*}
\begin{split}
& \left|\left( A^\epsilon \int_0^t e^{A(t-s)}BB^*P(\Phi(s)-I)x\, ds, z\right)_Y \right| 
\\[1mm]
& \qquad = \left|\int_0^t ( B^*P(\Phi(s)-I)x,B^*e^{A^*(t-s)}{A^*}^\epsilon z)_Y \,ds\right|
\\[1mm]
& \qquad \le  c \,\left[\int_0^t \|B^*P(\Phi(s)-I)x\|^{q'}_U\,ds\right]^{\frac1{q'}} 
\left[\int_0^t \|B^*e^{A^*\sigma}{A^*}^\epsilon z\|^q_U\, d\sigma\right]^{\frac1{q}}
\\[1mm]
& \qquad \le c \,\left[\int_0^t \|(\Phi(s)-I)x\|^{q'}_{\cD(A^\epsilon)}\,ds\right]^{\frac1{q'}} \|z\|_Y 
\\[1mm]
& \qquad \le c \,\left[\int_0^t s^{\alpha q'}\,ds \right]^{\frac1{q'}}\|z\|_Y 
\le c \,\|z\|_Y\, t^{\alpha+1/q'}\,.
\end{split}
\end{equation*}
Since $z$ is arbitrary, the claimed enhanced estimate follows. 

\medskip
Thus, still in the case $n_0>0$, we apply successively Claim~\ref{bootstrap} with 
\begin{equation*}
\alpha =\frac{1}{q'}, \, \frac{2}{q'},, \dots, \frac{n_0}{q'}\,,
\end{equation*}
until we obtain
\begin{equation}\label{ultima}
\left\|A^\epsilon \int_0^t e^{A(t-s)}BB^*P(\Phi(s)-I)x\, ds \right\|_Y 
\le c\, t^{(n_0+1)/q'}\,, \qquad 0\le t \le T\,.
\end{equation} 
Combining \eqref{ultima} with \eqref{stimabase} and noting that $n_0+1 > q'(1-\epsilon)$ by definition of $n_0$, we finally attain
\begin{equation}\label{stima_eps} 
\|(\Phi(t)-I)x\|_{D(A^\epsilon)} \le c \,\big(t^{(n_0+1)/q'}\vee t^{1-\epsilon}\big)
=c\,t^{1-\epsilon}\,, 
\quad 0\le t \le T\,.
\end{equation}
We note that the estimate \eqref{stima_eps} subsumes the case $n_0=0$ discussed separately in the beginning of this step.

\smallskip
\noindent
{\bf 4.}
We can now perform the final step. 
Suppose that $\epsilon<1/q'\wedge (1-\gamma)$. 
Then, a further use of Claim~\ref{bootstrap} with $\alpha = 1-\epsilon$ yields (from \eqref{stima_eps})
\begin{equation*}
\left\|A^\epsilon\int_0^t e^{A(t-s)}BB^*P(\Phi(s)-I)x\, ds\right\|_Y 
\le c\,t^{1-\epsilon + \frac1{q'}}\,, \qquad 0\le t \le T\,, 
\end{equation*}
where now $1-\epsilon + 1/q'>1$, owing to the assumed relation $\epsilon q'<1$. 
This implies, {\em a fortiori}, that
\begin{equation} \label{lim_zero}
\lim_{t\to 0^+} \frac1{t}\left\|\int_0^t e^{A(t-s)}BB^*P(\Phi(s)-I)x\, ds\right\|_Y =0\,.
\end{equation}

If $1/q'\le \epsilon < 1-\gamma$ instead, a more careful analysis is called for. 
For a given $z\in Y$ and any $t\in [0, T]$ we have
\begin{equation*}
\begin{split}
& \left(\int_0^t e^{A(t-s)}BB^*P(\Phi(s)-I)x\,ds, z\right)_Y
\\[1mm]
& \quad = \int_0^t \big(B^*P(\Phi(s)-I)x,B^*e^{A^*(t-s)}z\big)_Y\,ds 
\\[1mm]
& \quad = \int_0^t \big(B^*P(\Phi(s)-I)x,F(t-s)z\big)_Y\,ds 
+ \int_0^t \big(B^*P(\Phi(s)-I)x,G(t-s)z\big)_Y\,ds\,,
\end{split}
\end{equation*}
where in the last equality we have introduced the decomposition \eqref{e:key-hypo} of the adjoint of the kernel $e^{At}B$.
To estimate the first summand in the last row, we exploit the (singular, close to $0$) estimate \eqref{e:singular} pertaining to the component $F(t)$ of $B^*e^{A^*t}$, and recall once more the bound \eqref{stima_eps},
to find
\begin{equation} \label{e:bound_1}
\begin{split}
& \left|\int_0^t \left(B^*P(\Phi(s)-I)x,F(t-s)z\right)_Y\,ds\right|
\\[1mm]
& \qquad \le c \,\int_0^t \|(\Phi(s)-I)x\|_{\cD(A^\epsilon)} \|F(t-s)\|_{\cL(Y,U)} \,ds \, \|z\|_Y
\\[1mm] 
& \qquad \le c\,\int_0^t \frac{s^{1-\epsilon}}{(t-s)^\gamma} \,e^{-\eta(t-s)}\,ds \,\|z\|_Y
\le c\,t^{2-\epsilon -\gamma} \,\|z\|_Y\,.
\end{split}
\end{equation}
It is here that the constraint $\epsilon<1-\gamma$ arises, owing to the requirement $2-\epsilon -\gamma>1$.

As for the second summand, the basic regularity (in time) of the component $G(t)$ of $B^*e^{A^*t}$
(see the item 2.~of the Assumptions \ref{a:ipo_2}) combined with \eqref{stima_eps} brings about
\begin{equation} \label{e:bound_2}
\begin{split}
& \left|\int_0^t \left(B^*P(\Phi(s)-I)x,G(t-s)z\right)_Y\,ds\right|
\\[1mm]
& \qquad \le c \,\int_0^t \|(\Phi(s)-I)x\|_{\cD(A^\epsilon)} \|G(t-s)z\|_U \,ds 
\\[1mm]
& \qquad \le c \,\left[ \int_0^t s^{(1-\epsilon)p'} ds\right]^{1/p'} \|G(\cdot)\|_{\cL(Y,L^p(0,T,U))}\|z\|_Y 
\\
& \qquad \le c \,t^{1-\epsilon+1/p'} \|z\|_Y
\end{split}
\end{equation}
for any $t\in [0,T]$, and with the summability exponent $p\in [1,\infty)$.
Because $p\ge 1$ is finite and yet arbitrarily large, $1/p'$ can be taken arbitrarily close to $1$, 
which implies that $1-\epsilon+1/p'>1$. 
Since $z\in Y$ is also arbitrary, the estimates \eqref{e:bound_1} and \eqref{e:bound_2} yield
\begin{equation*}
\left\|\int_0^t e^{A(t-s)}BB^*P(\Phi(s)-I)x\, ds\right\|_Y \le c\,t^\alpha\,, \qquad 0\le t \le T\,,
\end{equation*}
with $\alpha>1$, and hence the limit \eqref{lim_zero} holds true, even in the case 
$1/q'\le \epsilon < 1-\gamma$.

\smallskip
\noindent
{\bf 5.}
We resume the representation \eqref{start} of the incremental ratio for the optimal state
semigroup, and combine the limit
\begin{equation*} 
\lim_{t\to 0^+} \frac{e^{At}-I}{t}(x-A^{-1}BB^*Px) = A(x-A^{-1}BB^*Px)  \quad \textrm{in $Y$}
\end{equation*}
(pointed out already in the step 1.) with \eqref{lim_zero}, to obtain that
\begin{equation*}
\lim_{t\to 0^+}\frac{\Phi(t)-I}{t}x = A(x-A^{-1}BB^*Px)\,,
\end{equation*}
whose meaning is $x\in \cD(A_P)$ and $A_Px= A(I-A^{-1}BB^*P)x$.
The set inclusion \eqref{e:seconda} is established, thus concluding the proof of Theorem~\ref{t:main}.

\qed


\section{The $\epsilon$-$\gamma$ constraint: PDE illustrations}
\label{s:illustrations}
We have seen in the previous section that we are able to characterize the domain of the optimal state semigroup's generator $A_P$, 
in the case the values of the (independent) parameters $\gamma$ and $\epsilon$
that occur in the standing Assumptions~\ref{a:ipo_2} fulfil the relation $\epsilon<1-\gamma$.
One may wonder whether this constraint is plausible, and actually met by significant boundary control systems which comprise hyperbolic and parabolic components and whose abstract formulation yields a kernel $e^{At}B$ (in fact, its adjoint) satisfying the Assumptions~\ref{a:ipo_2}.
Because the very requirement 3a) of the Assumptions~\ref{a:ipo_2} may be satisfied for a certain
$\epsilon_0>0$,
while being false (at least {\em a priori}) for $\epsilon<\epsilon_0$, the property $\epsilon<1-\gamma$ 
cannot be taken for granted.

In this section we show that the answer to the said question is positive. 
We recall and discuss briefly three distinct PDE problems describing distinct physical interactions:
from mechanical-thermal and acoustic-structure ones (in 2D and 3D, respectively) to fluid-elasticity ones in 3D.
Each of these illustrations has been previously shown to fall into the class of abstract control systems 
defined by the Assumptions~\ref{a:ipo_2}; see the papers \cite{abl-thermo_2005}, \cite{bucci_2008}, \cite{bucci-dampedfsi_2021}, respectively.
The values of $\gamma$ and $\epsilon$ specifically suited for either case had explicitly arisen during the proofs of the soughtafter (respective) trace regularity results that pertain to a certain component of the
solution to the uncontrolled PDE problem.

\medskip


\subsection{A thermoelastic system}
As our first example, we consider a classical PDE model which describes the displacements of an
elastic (thin) plate, in the presence of thermal effects; in regard to the modeling, see \cite{lagnese}.
This very system of coupled hyperbolic-parabolic PDE provided motivation for the distinctive abstract framework introduced in \cite{abl_2005}, along with the linear-quadratic theory devised therein,
subsequently continued in \cite{abl_2013,ab-uniqueness_2020} with focus on the infinite time horizon problem.

\smallskip
{\em The PDE problem.}
Let $\Omega\subset\mathbb{R}^2$ be a bounded domain with smooth boundary $\Gamma:=\partial\Omega$;
the symbol $\nu$ will denote the outward unit normal to the curve $\Gamma$.
The linear PDE model comprises a Kirchhoff equation for the vertical displacement $w(t,x)$ of the plate (and hence, $\rho>0$ below) and a simple diffusion equation for the temperature $\theta(t,x)$, with
the coupling occurring in the interior.
The system is supplemented with {\em clamped} (mechanical) boundary conditions and is subject to a {\em thermal} boundary control action $g$. 

Aiming to study the associated quadratic optimal control problem on an infinite time horizon,
we will have here $t\in (0,\infty)$.
Thus, the initial-boundary value problem (IBVP) is as follows:
\begin{equation} \label{e:thermo}
\begin{cases}
w_{tt} - \rho \Delta w_{tt} + \Delta^2 w + \Delta\theta=0 
&  \mbox{in } \,Q_\infty:=(0,\infty)\times \Omega
\\[1mm] 
\theta_t - \Delta \theta - \Delta w_t = 0
& \mbox{in } \,Q_\infty 
\\[1mm] 
w=\frac{\partial w}{\partial\nu}=0  
& \mbox{on } \Sigma_\infty:=(0,\infty)\times \Gamma
\\[1mm] 
\theta = g 
& \mbox{on } \Sigma_\infty 
\\[1mm]
w(0,\cdot) = w^0, \; w_t(0,\cdot) = w^1; \quad
\theta(0,\cdot) = \theta^0 &  \mbox{in } \Omega\,.
\end{cases}
\end{equation}

\medskip
{\em Functional setup, abstract formulation.}
The natural energy/state space is 
\begin{equation*}
Y=H^2_0(\Omega)\times H^1_0(\Omega)\times L^2(\Omega)\,.
\end{equation*}
With $U=L^2(\Gamma)$, we take $\cU:=L^2(0,\infty;U)$ as the space of admissible controls. 

The PDE problem \eqref{e:thermo} admits a reformulation as the Cauchy problem \eqref{e:cauchy-pbm},
after having set 
\begin{equation*}
y(t):=y(t,\cdot)\equiv (w(t,\cdot),w_t(t,\cdot),\theta(t,\cdot))\,,
\qquad 
g(t)=g(t,\cdot)\,,
\end{equation*}
and where $A$ and $B$ are suitable operators which can be detailed explicitly;
see \cite[Appendix~3J, p.~402]{las-trig-redbooks} and \cite{bucci-las-thermo_2004}.
%
We retrace below the various steps and analytical results which combine to bring about the conclusion that the PDE problem \eqref{e:thermo} falls within the abstract framework defined by the Assumptions~\ref{a:ipo_2}\footnote{Here we make reference to all the hypotheses with the only exception of 3b), that pertains to the observation operator $R$.}.

\smallskip
{\em Verification of the basic Assumptions~\ref{a:ipo_0}.}
(a) It is pretty well-known that the present (free) dynamics operator $A$ is the generator of a $C_0$-semigroup of contractions $e^{At}$ in $Y$, $t\ge 0$; this semigroup is {\em not} analytic, as opposed to the case when the thermoelastic system comprises the Euler-Bernoulli plate model ($\rho=0$), in place of the {\em hyperbolic} Kirchhoff one ($\rho>0$). 

(b) In addition, the semigroup $e^{At}$ is exponentially stable, as proved in \cite{avalos-las-thermo}
by energy methods.
(We note that the introduction of a novel (operatorial) multiplier is the decisive element in the proof therein.)

We recall that the study of well-posedness of initial-boundary value problems for (uncontrolled) linear thermoelastic systems has developed throughout the nineties of the last century, with initial focus on the case $k=0$ and {\em hinged} BC; see \cite[Chapter~3]{las-trig-redbooks} for an historical overview on the subject, along with a list of pertinent references.
It is worth noting that the analysis in \cite{las-trig-thermo_2000} offers a unified treatment of linear thermoelasticity under different BC, via an insightful decomposition of the underlying semigroup 
$e^{At}$ (which in particular covers the problem at hand here).

\smallskip
(c) As for the controlled thermoelastic problem \eqref{e:thermo}, the explicit computation of $B$ and 
$A^{-1}B$ performed in \cite{bucci-las-thermo_2004} -- combined with the aforementioned results -- confirms 
that the basic Assumption~\ref{a:ipo_0} are satisfied.

\smallskip
{\em Verification of the Assumptions~\ref{a:ipo_2}.} 
The applicability of the theory of the quadratic optimal control problem on an infinite time horizon 
devised in \cite{abl_2013} is discussed in full detail in \cite[\S\,2.1]{abl_2013}, as the end point
in a series of steps carried out in the earlier \cite{bucci-las-thermo_2004} and \cite{abl-thermo_2005}.
We recall specifically that the decomposition \eqref{e:key-hypo} of $B^*e^{A^*t}$, along with a
local (in time) version of the estimate \eqref{e:singular} for the component $F(t)$ 
(valid in a right neighbourhood of $t=0$, with $\eta=0$),
as well as the statements 2. and 3a) of the Assumptions~\ref{a:ipo_2} were shown to be valid 
in \cite{bucci-las-thermo_2004}.

The boundary regularity result that yields 3c) was subsequently sought and achieved in 
\cite[Theorem~1.1]{abl-thermo_2005}.
Finally, the enhanced estimate \eqref{e:singular} (for all $t>0$ and with $\eta>0$) was proved
in \cite[\S\,2.1]{abl_2013}.      
 
It is important to emphasize that the conditions 1., 2., 3a) on $B^*e^{A^*t}$ and 3c) on 
$B^*e^{A^*t}{A^*}^\epsilon$ have respective PDE counterparts; these are always regularity (in time and space) results for certain boundary traces\footnote{Below and later on we will utilize
the symbol $\leadsto$ to make explicit the said correspondence.}.
In the case of the thermoelastic system \eqref{e:thermo}, they specifically pertain to the normal traces (on $\Gamma$) of the thermal component of the {\em uncontrolled} system (i.e. \eqref{e:thermo} with $g\equiv 0$), as synthesized below:

\smallskip

\begin{itemize}

\item[$\bullet$]
regularity of the operator $B^*e^{A^*\cdot}$ $\leadsto$ 
regularity of $\frac{\partial \theta}{\partial\nu}\big|_{\Gamma}$,

\smallskip

\item[$\bullet$]
regularity of the operator $B^*e^{A^*\cdot}{A^*}^\epsilon$ $\leadsto$ 
regularity of $\frac{\partial \theta_t}{\partial\nu}\big|_{\Gamma}$.

\end{itemize}

\medskip
{\em Values of $\gamma$ and $\epsilon$.}
The outcome of \cite[Theorem~3.3]{bucci-las-thermo_2004} and \cite[Theorem~1.1]{abl-thermo_2005} allows to
conclude that the values of $\gamma$ and $\epsilon$ are as follows:
\begin{equation*}
\gamma= \frac{3}{4}+\sigma\,, \qquad 0<\epsilon <\frac{1}{4}\,, 
\end{equation*}
where $\sigma\in (0,1/4)$ may be arbitrarily close to $0$, while $\epsilon$ can be taken freely in the specified range.
Thus, since $1-\gamma\equiv 1/4-\sigma$, the constraint $\epsilon<1-\gamma$ is satisfied. 

\medskip


\subsection{A PDE model for acoustic-structure interactions with thermal effects} 
One may interpret the two dimensional plate dealt with in the previous paragraph 
as a {\em flat}, elastic portion $\Gamma_0$ of the boundary of a three dimensional, bounded domain $\Omega$ -- representing an acoustic chamber --, with the interactions between the acoustic waves in $\Omega$ and the vibrations of the flexible wall taking place on $\Gamma_0$.
Hence, differently than in the case of the thermoelastic system, 
the resulting structural acoustic model displays a {\em strong} coupling, as it occurs by means of boundary traces on $\Gamma_0$.
Optimal control problems arise naturally within this context, motivated by the goal of reducing the 
noise within the chamber and/or the vibrations of the elastic wall.

An archetypical PDE model for acoustic-structure interactions which couples equations of different type (a hyperbolic wave equation and a parabolic-like elastic equation) acting on manifolds of different dimensions, further subject to {\em point} control, has been the object of extensive research work during the nineties.
Point control mathematically describes -- via a combination of derivatives of Dirac delta functions -- the action of piezoelectric patches for the noise attenuation, as proposed by Dimitriadis~{\em et al.}
The work \cite{avalos-las-structural} by Avalos and Lasiecka studies the quadratic optimal control problem associated with the said linear problem.
It provides the first contribution to the mathematical analysis of Riccati equations (with unbounded operator coefficients) that arise from the optimal boundary control of systems of coupled hyperbolic-parabolic PDE
-- recast as abstract systems in the usual form $y'=Ay+Bg$ --, 
by exploiting (thus revealing) the central role of certain regularity estimates for the kernel $e^{At}B$, even in the {\em absence of analiticity} of the semigroup $e^{At}$.

Here we consider a PDE system which is a variant of the aforementioned structural acoustic model.
With the same basic equation for the acoustic waves propagation in the chamber,
the equation describing the vertical displacements of the elastic wall changes to embed rotational forces
(which account for the term $-\rho\Delta w_{tt}$ in \eqref{e:asi} below, $\rho>0$ being proportional to the thickness of the plate). 
And yet, the thermal effects -- which are present -- bring about an additional diffusion equation on
$\Gamma_0$. 
Thus, the resulting system comprises {\em three} evoutionary PDE, of different type: precisely, a hyperbolic equation and a system of coupled hyperbolic-parabolic PDE.

\medskip
{\em The PDE problem.}
Let $\Omega\subset \mathbb{R}^3$ be a bounded domain,  
with boundary 
\begin{equation*}
\partial \Omega=: \Gamma=\overline{\Gamma_0\cup \Gamma_1}\,, 
\end{equation*}
where $\Gamma_i\subset \mathbb{R}^2$, $i=1,2$, are open, simply connected and disjoint;
$\Gamma_1$ is the so called {\em hard} wall.
The PDE problem comprises a wave equation for the acoustic velocity potential
$z=z(t,x)$, $x\in \Omega$, 
and a thermoelastic system for the pair of the plate's vertical displacement and the temperature
$(w(t,x),\theta(t,x))$, $x\in \Gamma_0$.
The wave equation is supplemented with Neumann BC, while the thermoelastic system (is supplemented) with clamped BC and is subject to Dirichlet boundary control.
Thus, the IBVP is as follows:
\begin{equation} \label{e:asi} 
\begin{cases}
z_{tt} = \Delta z & \text{in $(0,\infty)\times \Omega=:Q_\infty$} 
\\[1mm]
\frac{\partial z}{\partial\tilde{\nu}} + d_1 z = 0 
& \text{on $(0,\infty)\times \Gamma_1=:\Sigma_1$}
\\[2mm]
\frac{\partial z}{\partial\tilde{\nu}} = w_t
& \text{on $(0,\infty)\times \Gamma_0=:\Sigma_0$}
\\[1mm]
w_{tt} - \rho \Delta w_{tt} + \Delta^2 w + \Delta\theta + z_t=0
& \text{in $\Sigma_0$}
\\[1mm] 
\theta_t - \Delta \theta - \Delta w_t = 0
&  \text{in $\Sigma_0$}
\\[2mm] 
w=\frac{\partial w}{\partial\nu}=0\,, \; \theta = g & \text{on $(0,\infty)\times \partial\Gamma_0$}
\\[2mm]
z(\cdot,0) = z^0, \; z_t(\cdot,0) = z^1 & \text{in $\Omega$}
\\[1mm]
w(\cdot,0) = w^0, \; w_t(\cdot,0) = w^1; \; \theta(\cdot,0) = \theta^0 & \text{in $\Gamma_0$.}
\end{cases}
\end{equation}
We note that symbols $\tilde{\nu}$ and $\nu$ above denote the outward unit normals to $\Gamma$ and to the curve $\partial\Gamma_0$, respectively;
$d_1$ and $\rho$ are positive constants ($\rho$ is proportional to the thickness of the plate).

\smallskip
{\em Functional setup, abstract formulation.}
We recall briefly that the natural (finite energy) state space for the PDE problem \eqref{e:asi} is 
\begin{equation*}
Y = H^1(\Omega)\times L^2(\Omega)\times H^2_0(\Gamma_0)\times H^1_0(\Gamma_0)\times L^2(\Gamma_0)\,;
\end{equation*}
with $U=L^2(\partial\Gamma_0)$, we take $\cU:=L^2(0,\infty;U)$ as the space of admissible controls. 
Indeed, the IBVP \eqref{e:asi} can be recast as the Cauchy problem \eqref{e:cauchy-pbm},
having set 
\begin{equation*}
y(t):=y(t,\cdot)\equiv (z(t,\cdot),w(t,\cdot),w_t(t,\cdot),\theta(t,\cdot))\,,
\qquad 
g(t)=g(t,\cdot)\,,
\end{equation*}
and with the operators $A$ and $B$ explicitly identified; see (2.13) and (2.15) in 
\cite[Section~2]{bucci_2007}.
Now, the complex of results -- be they functional-analytic or in the PDE realm -- showing that the control-theoretic properties listed as Assumptions~\ref{a:ipo_2} hold true (for the time being and as before, excluding the assumption 3b) on the observation operator $R$),  
are highlighted and properly attributed below.

\smallskip
{\em Verification of the basic Assumptions~\ref{a:ipo_0}.}
(a) 
The study of well-posedness of the abstract Cauchy problem 
\begin{equation*}
\begin{cases}
y'=Ay
\\[1mm]
y(0)=y_0\in Y
\end{cases}
\end{equation*}
corresponding to the IBVP \eqref{e:asi} with homogeneous boundary data (i.e., with $g\equiv 0$),
as well as of the uniform stability of the corresponding dynamics, has been pursued successfully by Lebiedzik in \cite{lebiedzik_2000}.
This work establishes that the said operator $A$ is the generator of a $C_0$-semigroup of contractions
$e^{At}$ in $Y$, $t\ge 0$. 
Classical semigroup theory provides the tool.
\\
(b) Furthermore, by using energy methods and a compactness-uniqueness argument, $e^{At}$ is shown to be exponentially stable.

\smallskip

(c) 
In the presence of nontrivial boundary data $g$, the control operator $B$ arises.
This and $A^{-1}B$ are explicitly computed by the second (named) author in \cite{bucci_2008}, thereby concluding that the basic Assumption~\ref{a:ipo_0} on $(A,B)$ are satisfied.

\smallskip
{\em Verification of the Assumptions~\ref{a:ipo_2}.}
The analysis performed in \cite{bucci_2008} pertains to the optimal control problem on a {\em finite} time horizon, meaning that $t$ varies on $[0,T)$ rather than on $[0,\infty)$; 
consequently, the asymptotic behaviour of $e^{At}$ -- as well as the one of the component $F(t)$ of 
$B^*e^{A^* t}$ -- does not play any role therein. 
However, \cite{bucci_2008} provides essentially the appropriate boundary regularity estimates that allow for the invocation of the theory in \cite{abl_2013}, besides the one in \cite{abl_2005} focused on the
case $T$ is finite.
We remind the reader that these pertain to the regularity (in time and space) of the boundary traces
of the thermal component, as summarized below: 

\smallskip
\begin{enumerate}

\item[$\bullet$]
regularity of the operator $B^*e^{A^*\cdot}$ $\leadsto$ 
regularity of $\frac{\partial \theta}{\partial\nu}\big|_{\partial\Gamma_0}$,

\vspace{1mm}

\item[$\bullet$]
regularity of the operator $B^*e^{A^*\cdot}{A^*}^\epsilon$ $\leadsto$ 
regularity of $\frac{\partial \theta_t}{\partial\nu}\big|_{\partial\Gamma_0}$.

\end{enumerate}
Precisely, Theorem~2.3 in \cite{bucci_2008} yields the sought decomposition \eqref{e:key-hypo} of
the operator $B^*e^{A^*t}$, along with the validity of 2., 3a) and 3c) of the Assumptions~\ref{a:ipo_2}.
A reworking on the proof of the singular estimate 
\begin{equation*}
\|F(t)\|_{\cL(Y,U)} \le C\,t^{-3/4-\sigma} \qquad \forall t>0
\end{equation*}
proved therein would allow to obtain the enhanced estimate \eqref{e:singular} (with appropriate constants $N$ and $\eta$), so that 1. of the Assumptions~\ref{a:ipo_2} holds true as well; the details are omitted.

\begin{remark}
\begin{rm}
We point out that the complex of (regularity) results provided by \cite[Theorem~2.3]{bucci_2008} follows combining

\begin{enumerate}

\item[-]
the sharp trace theory for the solutions to wave equations with nontrivial (Neumann) boundary data 
(in particular the estimate (3.9) in \cite{bucci_2008}, whose proof is found e.g. in \cite{las-trig-waves_1991}), which in turn yields

\item[-]
the improved boundary regularity of the elastic component of the (uncontrolled) system: precisely,
\begin{equation*}
\big\|\Delta w|_{\partial\Gamma_0}\big\|_{L^2(0,T;L^2(\partial\Gamma_0))}
\le C_T \|y_0\|_Y\,, \qquad y_0\in Y
\end{equation*}
(see \cite[Proposition~3.2]{bucci_2008}),

\item[-]
the smoothing effect of the parabolic component.

\end{enumerate}

\end{rm}
\end{remark}

\smallskip
{\em Values of $\gamma$ and $\epsilon$.}
In conclusion, the analysis pursued in \cite{bucci_2007,bucci_2008} and more specifically in
\cite[Theorem~2.3]{bucci_2008} bring about respective values of $\gamma$ and $\epsilon$ which are akin
to the ones found for the uncoupled thermoelastic system;
to wit,  
\begin{equation*}
\gamma= \frac{3}{4}+\sigma\,, \qquad 0<\epsilon <\frac{1}{4}\,, 
\end{equation*}
where $\sigma\in (0,1/4)$ may be arbitrarily close to $0$, while $\epsilon$ can be taken freely in the specified range.
Thus, since $1-\gamma\equiv 1/4-\sigma$, the constraint $\epsilon<1-\gamma$ is readily satisfied. 

\medskip


\subsection{A PDE model for a fluid-elasticity interaction} 
The PDE problem that serves as a last illustration of the feasibility of the relation
$\epsilon<1-\gamma$ (assumed in Theorem~\ref{t:main}) originates from the mathematical description of a very different physical interaction, more precisely a fluid-elasticity one.
The PDE problem under consideration arises from the linearization of a recognized PDE 
system that describes a fluid-structure interaction (FSI) in 3D, under the hypothesis of {\em stationary} interface.
The true nonlinear scenario is that of an elastic body fully immersed in a viscous, incompressible fluid;
the interface is the boundary of the elastic body, where the interactions between the two media
occur.
We refer the reader to the work of Barbu~{\em et al.}~\cite{barbu-etal_2007} and its bibliography, 
in regard to the original FSI.
As it will become apparent below, the said linearization further includes a dissipation term on the interface, which renders the coupled dynamics uniformly stable, as proved in \cite{avalos-trig-dampedfsi_2009}.

\smallskip
{\em The PDE problem.}
Let $\Omega\subset \mathbb{R}^3$ be the bounded and smooth domain representing the fluid-solid region.
If we denote by $\Omega_f$ and $\Omega_s$ the (open, smooth) domains occupied by the 
fluid and the solid, respectively, then $\Omega$ is the interior of 
$\overline{\Omega}_f\cup\overline{\Omega}_s$.
The velocity field of the fluid is represented by a vector-valued function $u$, which satisfies the
equations of Stokes flow in $\Omega_f$; the scalar function $p$ represents, as usual, the pressure.
The displacements of the solid region $\Omega_s$ are described by the variable $w$,
which satisfies the Lam\'e system of dynamic elasticity.
The interaction occurs at the interface, say, $\Gamma_s=\partial\Omega_s$, which is assumed stationary.
Then, $\Gamma_f$ will denote the outer boundary of $\Omega_f$ instead: namely, 
$\Gamma_f=\partial\Omega_f\setminus \partial\Omega_s$.
Below, $\nu=\nu(x)$ is the outward unit normal for the fluid region $\Omega_f$ and 
therefore it is pointing towards the interior of the solid region $\Omega_s$.
  
Thus, we consider the following system of PDE in the unknowns $(u,w)$:
\begin{equation} \label{e:dampedfsi}
\begin{cases} 
u_t-{\rm div}\,\epsilon(u)+\nabla p= 0 & \textrm{in $(0,\infty)\times \Omega_f=:Q_f$}
\\[1mm]
{\rm div}\, u=0 & \textrm{in $Q_f$}
\\[1mm] 
w_{tt} - {\rm div}\,\sigma(w)+w=0 &\textrm{in $(0,\infty)\times \Omega_s=:Q_s$}
\\[1mm]
u=0 & \textrm{on $(0,\infty)\times \Gamma_f=:\Sigma_f$}
\\[1mm]
\epsilon(u)\cdot\nu=\sigma(w)\cdot\nu + p\nu +g & \textrm{on $(0,\infty)\times \Gamma_s=:\Sigma_s$}
\\[1mm]
w_t-\sigma(w)\cdot\nu=u & \textrm{on $\Sigma_s$}
\\[1mm]
u(0,\cdot)=u_0 & \textrm{in $\Omega_f$}
\\[1mm]
w(0,\cdot)=w_0\,, \quad w_t(0,\cdot)=w_1 & \textrm{in $\Omega_s$.}
\end{cases}
\end{equation}
The symbols $\sigma$ and $\epsilon$ denote the elastic stress tensor and the strain tensor, respectively, that are
\begin{equation}\label{e:tensori}
\epsilon_{ij}(w)= \frac12\Big(\frac{\partial w_i}{\partial x_j}
+\frac{\partial w_j}{\partial x_i}\Big)\,,
\qquad
\sigma_{ij}(w) = \lambda \sum_{k=1}^3 \epsilon_{kk}(w)\delta_{ij}
+ 2\mu \epsilon_{ij}(w)\,,
\end{equation}
where $\lambda, \mu$ are the Lam\'e constants and $\delta_{ij}$ is the Kronecker symbol.
(We note the abuse of notation determined by the use of the letter $\epsilon$ for the parameter in 3. of
the Assumptions~\ref{a:ipo_2}, and also for the elastic strain tensor defined in \eqref{e:tensori}.
However, as the latter appears only in \eqref{e:dampedfsi}, no confusion is likely to arise.) 

The PDE system in \eqref{e:dampedfsi} shares with the previous illustrations the features that it comprises both parabolic and hyperbolic equations, and moreover the parabolic component is subjected to a boundary control action.
Notably, its abstract reformulation falls into the framework characterized by the Assumptions~\ref{a:ipo_2},
as proved in \cite{bucci-dampedfsi_2021}.

\smallskip
{\em Functional setup, abstract formulation.}
The functional setup for the IBVP \eqref{e:dampedfsi} follows -- in its basic elements -- the
one introduced in \cite{barbu-etal_2007} for the true nonlinear problem.
Also, the perspective taken in \cite{bucci-dampedfsi_2021} to attain the sought regularity results
is akin to the one adopted in \cite{las-tuff-fsi_2009}, \cite{bucci-las-fsi_2010} 
and \cite{bucci-las-fsi_2011} in the study of the optimal boundary control problem associated with its (undamped) linearization.
We note that the regularity analysis that has been carried out in the aforementioned works 
benefited from the enhanced regularity of the boundary traces of the normal component of the stress tensor
on the interface, valid for weak as well as semigroup solutions.

The state (energy) space for the system is 
\begin{equation*}
Y := \cH\times [H^1(\Omega_s)]^d\times [L^2(\Omega_s)]^d\,,
\end{equation*}

\noindent
where $\cH$ denotes the ({\em null-div}) space pertaining to the fluid component, that is
\begin{equation*}
\cH := \big\{ u\in [L^2(\Omega_f)]^d\colon \; 
{\rm div}\, u=0\,, \; u\cdot \nu|_{\Gamma_f}=0\big\}\,,
\end{equation*}
while $[H^1(\Omega_s)]^3\times [L^2(\Omega_s)]^3$ is the natural energy space for
the Lam\'e system. 

\smallskip
As before, the IBVP \eqref{e:dampedfsi} can be recast as the abstract Cauchy problem \eqref{e:cauchy-pbm},
having set 
\begin{equation*}
y(t):=y(t,\cdot)\equiv (u(t,\cdot),w(t,\cdot),w_t(t,\cdot))\,,
\qquad 
g(t)=g(t,\cdot)\,,
\end{equation*}
and with suitable operators $A$ and $B$ explicitly identified; see \cite[Section~2.1]{bucci-dampedfsi_2021}.
The abstract setup introduced in \cite{bucci-dampedfsi_2021} has proved well-suited for the 
derivation of the boundary regularity estimates that allow for the invocation of
the theory of the infinite time horizon optimal control problem devised in \cite{abl_2013}.

\smallskip
{\em Verification of the basic Assumptions~\ref{a:ipo_0}.}
The basic Assumptions~\ref{a:ipo_0} for the pair of the (free) dynamics operator $A$ and the 
control operator $B$ are ascertained in \cite[Proposition~2.2]{bucci-dampedfsi_2021}.
The modified transmission condition $$w_t-\sigma(w)\cdot\nu=u$$ on the interface $\Gamma_s$ -- in place of 
$w_t=u$, which conveys the matching of the two velocities -- subsumes a source of dissipation, 
which renders the PDE system uniformly (exponentially) stable, as proved in 
\cite{avalos-trig-dampedfsi_2009}.

\smallskip
{\em Verification of the Assumptions~\ref{a:ipo_2}.}
The requirements 1., 2., 3a) on $B^*e^{A^*t}$ and 3c) on $B^*e^{A^*t}{A^*}^\epsilon$ of the Assumptions~\ref{a:ipo_2} translate into appropriate regularity results (in time and space) for 
the boundary traces of the fluid velocity field $u$ on $\Gamma_s$.
More specifically,

\smallskip
\begin{enumerate}

\item[$\bullet$]
regularity of the operator $B^*e^{A^*\cdot}$ $\leadsto$ 
regularity of $\frac{\partial u}{\partial\nu}\big|_{\Gamma_s}$,

\vspace{1mm}

\item[$\bullet$]
regularity of the operator $B^*e^{A^*\cdot}{A^*}^\epsilon$ $\leadsto$
regularity of $\frac{\partial u_t}{\partial\nu}\big|_{\Gamma_s}$,

\end{enumerate}
where $(u,w,w_t)$ solves the {\em uncontrolled} problem, that is IBVP \eqref{e:dampedfsi} with $g\equiv 0$.

\smallskip
It was found in \cite{bucci-dampedfsi_2021} that the sought boundary regularity results for the fluid
velocity field can be achieved exploiting and combining carefully the parabolic regularity of the fluid variable with the exceptional regularity of certain normal traces of the elastic component.
The latter, enhanced by the presence of the boundary dissipation into the coupled system,
emerges from an energy equality obtained in \cite{avalos-trig-dampedfsi_2009}, which in turn can be proved with the only use of multiplier methods.
In this respect it should be noted that the prior work carried out in \cite{las-tuff-fsi_2009} and 
\cite{bucci-las-fsi_2010,bucci-las-fsi_2011} on the (undamped) linearizations of the original FSI necessitated arguments from {\em microlocal analysis}, instead; see \cite[Lemmas~2.3 and 2.4]{bucci-las-fsi_2011}.

\smallskip
{\em Values of $\gamma$ and $\epsilon$.}
The analysis pursued in \cite{bucci-dampedfsi_2021} (whose core is the proof of Theorem~1.2 given in Section~2.3 therein), yields the following values of $\gamma$ and $\epsilon$:
\begin{equation*}
\gamma= \frac{1}{4}+\sigma\,, \qquad 0<\epsilon <\frac{1}{2}\,, 
\end{equation*}
where $\sigma\in (0,3/4)$ may be arbitrarily close to $0$, while $\epsilon$ can be taken freely in the specified range.
Thus, since $1-\gamma\equiv 3/4-\sigma$, the constraint $\epsilon<1-\gamma$ is readily satisfied. 


{\small
\section*{Acknowledgements}
\noindent
The research of P.~Acquistapace was partially supported by the PRIN-MIUR Project 2017FKHBA8 of the
Italian Education, University and Research Ministry.
\\
The research of F.~Bucci was partially supported by the Universit\`a degli Studi di Firenze
under the 2020 Project {\em Metodi ed Applicazioni per Equazioni Differenziali Ordinarie e a Derivate Parziali}.
Bucci is a member of the Gruppo Nazionale per l'Analisi Mate\-ma\-tica, la Probabilit\`a 
e le loro Applicazioni (GNAMPA) of the Istituto Nazionale di Alta Matematica (INdAM),
and participant to the GNAMPA Project 
{\em Problemi inversi e di controllo per equazioni di evoluzione e loro applicazioni} (2020).
She is also a member of the French-German-Italian Laboratoire International
Associ\'e (LIA) COPDESC in Applied Analysis.
}

\medskip


\end{document}